\documentstyle[12pt]{article}
\def\Bbb R{{\rm \bf R}}
\def\proclaim#1{\vskip2mm{\bf #1}\em}
\def\endproclaim{\em \vskip2mm}
\def\tag#1{\eqno(#1)}
\def\gathered{\begin{array}{c}}
\def\endgathered{\end{array}}
\def\text{\mbox}

\begin{document}

\title {An inverse source problem for the heat equation\\
and\\
the enclosure method}
\author{Masaru IKEHATA\\
Department of Mathematics,
Faculty of Engineering\\
Gunma University, Kiryu 376-8515, JAPAN}
\date{2 February 2006}
\maketitle
\begin{abstract}
An inverse source problem for the heat equation is considered.
Extraction formulae for information about the time and location
when and where the unknown source of the equation firstly appeared
are given from a single lateral boundary measurement. New roles of
the plane progressive wave solutions or their complex versions for
the backward heat equation are given.

\noindent
AMS: 35R30, 80A23

\end{abstract}

\section{Introduction}

Let $\Omega$ be a bounded domain of $\Bbb R^n (n=1,2,3)$ with smooth boundary.
Let $T$ be an arbitrary positive number. Let $u=u(x,t)$ satisfy
$$\displaystyle
u_t=\triangle u+f(x,t)\,\,\text{in}\,\Omega\times]0,\,T[.
$$

In this paper we consider an inverse problem for the heat equation.
The problem is

{\bf\noindent Inverse Problem.}
Assume that there exist a non negative number $T_0$ less than $T$ and
point $x_0\in\Omega$ such that $f(x_0,T_0)\not=0$ and
$f(x,t)=0$ for all $0<t<T_0$ and all $x\in\Omega$.

\noindent
{\it Extract} $T_0$ and information about the set $\{x\,\in\Omega\,\vert\,f(x,T_0)\not=0\}$
from the {\it data} $u\vert_{\partial\Omega\times]0,\,T[}$,
$\partial u/\partial\nu\vert_{\partial\Omega\times]0,\,T[}$
and $u(\,\cdot\,,0)$.

\noindent The number $T_0$ and the set
$\{x\,\in\Omega\,\vert\,f(x,T_0)\not=0\}$ are the time and
position when and where the heat source $f(x,t)$ firstly appeared.
Note that one may assume that $u(\,\cdot\,,0)=0$ in $\Omega$ and
in what follows we do so.

\noindent
Some related known results should be mentioned.

Yamatani-Ohnaka \cite{YO} considered the case when the source $f(x,t)$ takes the form
$$\displaystyle
f(x,t)=\sum_{j=1}^N p_j\delta(x-x_j,t-t_j)
$$
where the points $x_j\in\,\Omega$, times $t_j\in\,]0,\,T[$, strength $p_j(<0)$ and $N$ are all {\it unknown}.

\noindent
Note that, in this case
$T_0=\min_j\{t_j\}$, formally the set $\{x\in\Omega\,\vert f(x,T_0)\not=0\}$ coincides with $\{x_j\,\vert\,T_0=t_j\}$
and the source becomes {\it inactive} after the time $\max_j\,\,{t_j\,}$.
They made use of special solutions of the {\it backward heat equation} $v_t+\triangle v=0$
in $\Omega\times\,]0,\,T'[$ where $T\ge T'>\max_j\,{t_j\,}$ having a large parameter $c$.
Those solutions have singularity
at $t=T'$ on the given plane outside $\overline\Omega$
and are vanishing at $t\uparrow T'$ on $\overline\Omega$.
Those are constructed from the fundamental solution of the backward heat equation
in one-dimensional space.  Then, using integration by parts, they obtain a system of
equations involving those solutions, unknown $x_j, t_j, p_j$ and $N$. Carefully analyzing the asymptotic behaviour
of the system as $c\longrightarrow\infty$,
they gave a reconstruction formula of the source itself.

Yamamoto \cite{Y} considered an inverse source problem for
the wave equation, however, it is possible to apply his method to
the inverse source problem for the heat equation.
One interpretation of his method for the heat equation is the following.
The assumption on the heat source is that
$$\displaystyle
f(x,t)=\sigma(t)f(x)
$$
where $\sigma$ is a {\it known} function with $\sigma(0)\not=0$ and $f$ is {\it unknown}.
From our point of view this is the case when $T_0$ is {\it known} and $T_0=0$.

\noindent
Then, his method consists of two parts:

$\bullet$  a way of calculating the data
$\partial w/\partial\nu\vert_{\partial\Omega}\times\,]0,\,T[$
where $w$ solves
$$\begin{array}{c}
\displaystyle
w_{t}=\triangle w\,\,\text{in}\,\Omega\times\,]0,\,T[,\\
\\
\displaystyle
w(x,0)=f(x)\,\,
\text{in}\,\Omega,\\
\\
\displaystyle
w(x,t)=0\,\,\text{on}\,\partial\Omega\times\,]0,\,T[,
\end{array}
$$
from the data $u\vert_{\partial\Omega}\times\,]0,\,T[$,
$\partial
u/\partial\nu\vert_{\partial\Omega}\times\,]0,\,T[$
and the {\it time derivative}
$(
\partial/\partial t)
\partial u/\partial
\nu\vert_{\partial\Omega}\times\,]0,\,T[$;

$\bullet$  {\it reconstruction formula} of $f(x)$ itself from the calculated data
$\partial w/\partial\nu\vert_{\partial\Omega}\times\,]0,\,T[$.

\noindent
The first part can be done by solving a Volterra equation of the second kind.
For establishing the second part the {\it null-controllability} for the heat
equation \cite{FI, LR} and the completeness of the eigenfunctions of
the Dirichlet Laplacian in $\Omega$ are essential.

In \cite{EH2} El Badia-Ha Duong also considered the case when
$T_0=0$. They assume that: the source is supported
in a subset of $\Omega$ (mainly the point sources) that is independent of the time,
the strength of the source on the set
does not depend on the position and has a definite signature in an interval
$]0,\,T_{\star}[$ with $T_{\star}<T$; the source becomes {\it
inactive} after the time $T_{\star}$.

\noindent
Their method also consists of two parts:

$\bullet$  a way of calculating the temperature in
$\Omega$ at $t=T_1$ from the data
$u\vert_{\partial\Omega}\times\,]T_1,\,T[$ and $\partial
u/\partial\nu\vert_{\partial\Omega}\times\,]T_1,\,T[$ where $T_1$
is an arbitrary fixed time with $T_{\star}<T_1<T$;

$\bullet$  extracting moments of the source from $u(x,T_1)$ and the data
and determination of the sources form the moments.

\noindent The first part is also an application of the
null-controllability for the heat equation and the completeness of
the eigenfunctions of the Dirichlet Laplacian in $\Omega$.
\noindent
In the second part they start with the determination
of the number and location of point sources.
This is purely algebraic and can be done by using a previous result in the
inverse point source problem
for elliptic equation \cite{EH1}.
The next step is the determination of the strength of the point sources
from some exponential moments which is a combination of the
M\"untz's theorem (e.g., \cite{GL}) and an invertibility of a
matrix involving special solutions of the Helmholtz equations. As
they pointed out the choice of the matrix is not constructive.
These steps can be done by a combination of
an integration by parts and special solutions
for the backward heat equation.

It should be pointed out that
the controllability for the heat equation plays an important role also
in the problem of recovering the conductivity coefficient of the heat equation,
see Avdonin-Belishev-Rozkov \cite{ABR}.

In this paper we apply the {\it enclosure method} to {\bf Inverse Problem}
which was introduced by Ikehata \cite{I4} for inverse boundary value problems
for the elliptic equations.
In \cite{I1} Ikehata considered in two-dimensions, an
inverse source problem for the inhomogeneous Helmholtz equation $\triangle u+k^2u=f$
with the unknown source $f$ having the form $\chi_D(x)\rho(x)$ and
established an extraction formula of the convex
hull of $D$ provided {\it both} $D$ and $\rho$ are {\it unknown} and
$D$ is {\it polygonal}.
Now the enclosure method has been applied to inverse boundary
value problems \cite{I2,I3,IO}, a Cauchy problem for the stationary
Schr\"odinger equation \cite{I5} and inverse obstacle scattering problems
\cite{I6, I7, I8}. The common character of these problems is: the
governing equation is elliptic.

It is quite curious to consider whether one can apply the idea of the
enclosure method to the case when the governing equation is {\it
non elliptic} since there are many inverse problems whose governing equations are
non elliptic.  The inverse source problem for the heat equation is a typical and important one.
The aim of this paper is to introduce a direct approach
to the inverse source problem for the heat equation by employing the idea of the enclosure
method.

\noindent
The main feature of our method is:

$\bullet$  it is based on a simple {\it one line} formula;

$\bullet$  we do not make use of the exact controllability of the heat
equation nor the completeness of the eigenfunctions of
the Dirichlet Laplacian in $\Omega$;

$\bullet$  the assumptions on the unknown source is quite general;

$\bullet$  the method provides us a brief information about the
time and the location when and where the unknown source firstly
appeared instead of the detailed information of the source.

\noindent
For several other formulations of the inverse source
problem for the heat equation and results including uniqueness and
stability see Cannon-Esteva \cite{CE}, Cannon-DuChateau \cite{CD},
Isakov \cite{Is}, Hetllich-Rundell \cite{HR}, recent
publication of Trong-Long-Alain \cite{TLA} and references therein.

\section{Statement of results}

The results of this paper are divided into two parts. The first
one is a formula for extracting $T_0$ only and second is for
extracting both $T_0$ and information about the location of the set
$\{x\in\,\Omega\,\vert f(x,T_0)\not=0\}$.

{\bf\noindent 2.1.  Extracting $T_0$.}

The list of the assumptions on the source $f(x,t)$ is
the following.

(A.1)  $f(x,t)$ takes the form
$$\displaystyle
f(x,t)=\chi_D(x,t)\rho(x,t)
$$
where $D\subset\overline\Omega\times\,[T_0,\,T]$ is a Lebesgue measurable set;
$\rho$ is essentially bounded on $D$.

\noindent
Thus the source may appear at the time $T_0$ firstly and the points in $D\cap(\Bbb R^n\times\{T_0\})$.

\noindent
The next assumption describes the standing behaviour of the source.

(A.2)  There exist positive numbers
$\delta<T-T_0$, $C_1$, $C_2$ and $p\in\,[0,\,\infty[$
such that the $n$-dimensional Lebesgue
measure of the set $D(s)\equiv D\cap(\Bbb R^n\times\{T_0+s\})$
has the lower estimate:
$$\displaystyle
\vert D(s)\vert\ge C_1s^p\,\,\text{for almost all}\, s\in\,[0,\,\delta_1];
$$
the strength of the source $\rho$ satisfies
$$\displaystyle
\rho(x,t)\ge C_2\,\,\text{for almost all}\,(x,t)\in\,D\cap(\Bbb R^n\times\,[T_0,\,T_0+\delta])
$$
or
$$\displaystyle
-\rho(x,t)\ge C_2\,\,\text{for almost all}\,(x,t)\in\,D\cap(\Bbb R^n\times\,[T_0,\,T_0+\delta]).
$$

{\bf\noindent Definition 2.1.}
Given $\omega\in S^{n-1}$, $s\in\Bbb R$ define the {\it indicator function}
$I_{\omega}(\tau;s)$ by the formula
$$\displaystyle
I_{\omega}(\tau;s)
=e^{\tau s}\int_0^{T}\int_{\partial\Omega}\left(\frac{\partial v}{\partial\nu}u
-\frac{\partial u}{\partial\nu}v\right)dSdt,\,\,\tau>0
$$
where $v=v(x,t)=e^{\sqrt{\tau}x\cdot\omega-\tau t}$.

\noindent
The function $v$ takes positive values and satisfies the {\it backward} heat equation $v_t+\triangle v=0$
in the whole space. Moreover $e^{\tau s}v$ has the special character:

$\bullet$ if $t>s$, then $\lim_{\tau\longrightarrow\infty}e^{\tau s}v(x,t)=0$;

$\bullet$ if $t<s$, then $\lim_{\tau\longrightarrow\infty}e^{\tau s}v(x,t)=\infty$.

\noindent
Note that the function $w(x)=e^{\sqrt{\tau}x\cdot\omega}$
satisfies the equation
$$\displaystyle
\triangle w-\tau w=0
$$
in the whole space.  Since $\tau>0$, this is not the Helmholtz equation.
In the following we make use of the asymptotic behaviour of $w$ as
$\tau\longrightarrow\infty$.

\proclaim{\noindent Theorem 2.1.}
Assume that assumptions (A.1), (A.2) are satisfied.
As $\tau\longrightarrow\infty$ we have
$$\begin{array}{c}
\displaystyle
\lim_{\tau\longrightarrow\infty}\frac{\log\vert I_{\omega}(\tau;0)\vert}
{\tau}
=-T_0.
\end{array}
\tag {2.1}
$$
\noindent
Moreover, we have:

if $s<T_0$, then $\lim_{\tau\longrightarrow\infty}\vert
I_{\omega}(\tau;s)\vert=0$;

if $s>T_0$, then $\lim_{\tau\longrightarrow\infty}\vert
I_{\omega}(\tau;s)\vert=\infty$.

\endproclaim

The result can be generalized to the case when the governing equation
has {\it variable coefficients}.  In that case there is an interesting relationship between
the high frequency {\it asymptotic solution} of the corresponding wave equation \cite{ES} and that of
the formal adjoint equation of the heat equation.  This is clarified
in Section 3.

{\bf\noindent 2.2.  Extracting both $T_0$ and information about the set $\{x\,\in\Omega\,\vert\,f(x,T_0)\not=0\}$.}

This subsection starts with the following observation.

\noindent Let $n\ge 1$.  Let $c\not=0$ be an arbitrary number. Let
$$\displaystyle
v(x,t)=e^{-(z\cdot z)t}e^{x\cdot z}
\tag {2.2}
$$
where $\tau$ satisfies $\tau>c^{-2}$ and
$$
z=\left\{\begin{array}{lr}
\displaystyle
c\tau\left(\omega+i\sqrt{1-\frac{1}{c^2\tau}}\omega^{\perp}\right), & \quad n\ge 2,\\
\\
\displaystyle
c\tau\left(1+i\sqrt{1-\frac{1}{c^2\tau}}\,\right), & \quad n=1.
\end{array}
\right.
$$

\noindent
The function $v$ is a complex valued function and satisfies the backward heat equation $v_t+\triangle v=0$.
Moreover $e^{\tau s}v$ has the special character:

$\bullet$ if $s<t-cx\cdot\omega$, then $\lim_{\tau\longrightarrow\infty}e^{\tau s}\vert v(x,t)\vert=0$;

$\bullet$ if $s>t-cx\cdot\omega$, then $\lim_{\tau\longrightarrow\infty}e^{\tau s}\vert v(x,t)\vert=\infty$.

\noindent So one can expect that: using this function, one may obtain
more information about the location and shape of $D$ in the space
time. We show that the method \cite{I5} is applicable
in the case when $n=1,2$ (at least).

{\bf\noindent Definition 2.2.}  Let $n\ge 2$.
Given $c>0$, $s\in\Bbb R$, $\omega,\,\omega^{\perp}\in S^{n-1}$ with $\omega\cdot\omega^{\perp}=0$
define the {\it indicator function}
$I_{\omega,\,\omega^{\perp},c}(\tau;s)$ by the formula
$$\displaystyle
I_{\omega,\,\omega^{\perp},c}(\tau;s)
=e^{\tau s}\int_0^{T}\int_{\partial\Omega}\left(\frac{\partial v}{\partial\nu}u
-\frac{\partial u}{\partial\nu}v\right)dSdt,\,\,\tau>c^{-2}
$$
where $v$ is the function given by (2.2).

\noindent
Let $n=1$.
Given $c\not=0$, $s\in\Bbb R$,
define the {\it indicator function}
$I_{c}(\tau;s)$ by the formula
$$\displaystyle
I_{c}(\tau;s)
=e^{\tau s}\int_0^{T}\int_{\partial\Omega}\left(\frac{\partial v}{\partial\nu}u
-\frac{\partial u}{\partial\nu}v\right)dSdt,\,\,\tau>c^{-2}
$$
where $v$ is the function given by (2.2).

\noindent

Let $n\ge 2$.
Given $\omega\in\,S^{n-1}$ and $c>0$ define
the unit vector $\omega(c)$ in $\Bbb R^{n+1}=\Bbb R^n\times\Bbb R$ directed
into the half space $t<0$ by the formula
$$\displaystyle
\omega(c)=\frac{1}{\sqrt{c^2+1}}
\left(\begin{array}{c}
\displaystyle c\omega\\
\displaystyle
-1
\end{array}
\right).
$$
Then
$$
\displaystyle
t-cx\cdot\omega=-\sqrt{c^2+1}\left(\begin{array}{c} x\\ t\end{array}\right)\cdot\omega(c).
$$
Thus we see that:

$\bullet$  if $-s/\sqrt{c^2+1}>(x\,\,t)^T\cdot\omega(c)$, then
$\lim_{\tau\longrightarrow\infty}e^{\tau s}\vert v(x,t)\vert=0$;

$\bullet$  if $-s/\sqrt{c^2+1}<(x\,\,t)^T\cdot\omega(c)$, then
then $\lim_{\tau\longrightarrow\infty}e^{\tau s}\vert v(x,t)\vert=\infty$;

\noindent
First we consider the case when $n=2$.

\noindent
We assume that the unknown source takes the form
$$\displaystyle
f(x,t)=\sum_{j=1}^N\chi_{P_j\times\,[T_j,\,T]}(x,t)\rho_j(x,t)
$$
where

$\bullet$  $P_j\subset\overline\Omega$ is given by the interior of a
{\it polygon} and $T_j$ satisfies $0\le T_j<T$;

$\bullet$  if $j\not=j'$, then $\overline P_j\cap\overline P_{j'}=\emptyset$;

$\bullet$  for each $j$ $\rho_j\in L^{\infty}(P_j\times]T_j,\,T[)$ and for each vertices $p$ of the convex hull of
$P_j$ $\rho_j$ coincides with a H\"older
continuous function with exponent $\theta_j\in\,]0,\,1]$ in a neighbourhood of $(p,T_j)$
that does not vanish at $(p,T_j)$.

\noindent
We set $\displaystyle D=\cup_{j=1}^N(P_j\times\,]T_j,\,T[)$.

\proclaim{\noindent Theorem 2.2.}
Let $\omega(c)$ be regular with respect to $D$.
Assume that
$$\displaystyle
\sup_{x\in\,\Omega}\left(\begin{array}{c}
x\\
T\end{array} \right) \cdot\omega(c)<h_D(\omega(c)).
\tag {2.3}
$$

\noindent
Then the formula
$$\displaystyle
\lim_{\tau\longrightarrow\infty}
\frac{\displaystyle\log\vert I_{\omega,\,\omega^{\perp},c}(\tau;0)\vert}
{\tau}
=\sqrt{c^2+1}\,h_D(\omega(c)),
$$
is valid.

\noindent
Moreover, we have:

if $h_D(\omega(c))\le -s/\sqrt{c^2+1}$, then
$\lim_{\tau\longrightarrow\infty}\vert
I_{\omega,\,\omega^{\perp},c}(\tau;s)\vert=0$;

if $h_D(\omega(c))>-s/\sqrt{c^2+1}$, then
$\lim_{\tau\longrightarrow\infty}\vert
I_{\omega,\,\omega^{\perp},c}(\tau;s)\vert=\infty$.

\endproclaim

\noindent The condition (2.3) means that the set
$\Omega\times\{T\}$ is located in the half space
$(x,t)^T\cdot\omega(c)<h_D(\omega(c))$.  This can be satisfied in
the case when: $T$ is given and $c$ is sufficiently small; $c$ is
given and $T$ is sufficiently large.

Next we consider the case when $n=1$.

\noindent
We assume that the unknown source takes the form
$$\displaystyle
f(x,t)=\sum_{j=1}^N\chi_{P_j\times\,[T_j,\,T]}(x,t)\rho_j(x,t)
$$
where

$\bullet$  $P_j\subset\overline\Omega$ is given by the interior of an
{\it interval} and $T_j$ satisfies $0\le T_j<T$;

$\bullet$  if $j\not=j'$, then $\overline P_j\cap\overline P_{j'}=\emptyset$;

$\bullet$  for each $j$ $\rho_j\in L^{\infty}(P_j\times]T_j,T[)$ and for each end point $p$ of $P_j$
$\rho_j$ coincides with a function of class $C^2$
in a neighbourhood of the point $(p,T_j)$ that does not vanish at $(p, T_j)$.

\noindent
We set $\displaystyle D=\cup_{j=1}^N(P_j\times\,]T_j,\,T[)$.

Given $c\not=0$ define
the unit vector $\omega(c)$ in $\Bbb R^{1+1}=\Bbb R\times\Bbb R$ directed
into the half space $t<0$ by the formula
$$\displaystyle
\omega(c)=\frac{1}{\sqrt{c^2+1}}
\left(\begin{array}{c}
\displaystyle c\\
\displaystyle
-1
\end{array}
\right).
$$

\proclaim{\noindent Theorem 2.3.}
Let $\omega(c)$ be regular with respect to $D$.
Assume that
$$\displaystyle
\sup_{x\in\,\Omega}\left(\begin{array}{c}
x\\
T\end{array} \right) \cdot\omega(c)<h_D(\omega(c)).
$$

\noindent
Then the formula
$$\displaystyle
\lim_{\tau\longrightarrow\infty}
\frac{\displaystyle\log\vert I_{c}(\tau;0)\vert}
{\tau}
=\sqrt{c^2+1}\,h_D(\omega(c)),
$$
is valid.

\noindent
Moreover, we have:

if $h_D(\omega(c))\le -s/\sqrt{c^2+1}$, then
$\lim_{\tau\longrightarrow\infty}\vert I_{c}(\tau;s)\vert=0$;

if $h_D(\omega(c))>-s/\sqrt{c^2+1}$, then
$\lim_{\tau\longrightarrow\infty}\vert I_{c}(\tau;s)\vert=\infty$.

\endproclaim

\noindent
If $\Omega\subset\Bbb R^3$, the situation seems more
complicated.
We will consider the case in the future work.

\section{Proof of Theorem 2.1 and its generalization}

{\bf\noindent 3.1.  Proof of Theorem 2.1.}

Since we have the trivial identity
$$\displaystyle
I_{\omega}(\tau;s)=e^{\tau(s-T_0)}I_{\omega}(\tau;T_0)\,\,\forall\,s,
\tag {3.1}
$$
it suffices to study the asymptotic behaviour of the indicator function
at $s=T_0$.
Integration by parts gives
$$\displaystyle
\int_0^T\int_{\Omega}f(x,t)v(x,t)dxdt
=\int_{\Omega}u(x,T)v(x,T)dx
+\int_0^T\int_{\partial\Omega}
(\frac{\partial v}{\partial\nu}u-\frac{\partial u}{\partial\nu}v)dSdt.
$$
Thus we have the representation of the indicator function at $s=T_0$:
$$\displaystyle
I_{\omega}(\tau;T_0)
=e^{\tau T_0}\int_0^T\int_{\Omega}f(x,t)v(x,t)dxdt
-e^{\tau T_0}\int_{\Omega}u(x,T)v(x,T)dx.
\tag {3.2}
$$
Since $T>T_0$, it is easy to see that the absolute value of the second term of the right hand side
is dominated by
$$\displaystyle
\int_{\Omega}\vert u(x,T)\vert dxe^{\sqrt{\tau}\sup_{x\in\Omega}x\cdot\omega}e^{-\tau(T-T_0)}
=O(e^{-\tau(T-T_0)/2})
\tag {3.3}
$$
as $\tau\longrightarrow\infty$.

\noindent
On the other hand, from the assumptions (A.1) and (A.2)
one knows that the first term of the right hand side of (3.2) takes the form
$$
\begin{array}{c}
\displaystyle
e^{\tau T_0}\int_{D}\rho(x,t)v(x,t)dxdt\\
\\
\displaystyle
=e^{\tau T_0}\int_{D\cap\,(\Bbb R^n\times\,[T_0,\,T_0+\delta])}\rho(x,t)e^{\sqrt{\tau}x\cdot\omega-\tau t}dxdt
+e^{\tau T_0}\int_{D\setminus
(\Bbb R^n\times\,[T_0,\,T_0+\delta])
}\rho(x,t)e^{\sqrt{\tau}x\cdot\omega-\tau t}dxdt.
\end{array}
\tag {3.4}
$$
Since the set $D\setminus(\Bbb R^n\times\,[T_0,\,T_0+\delta]$ is contained in $\Bbb R^n\times
\,[T_0+\delta,\,T]$, we have
$$\displaystyle
e^{\tau T_0}\int_{D\setminus
(\Bbb R^n\times\,[T_0,\,T_0+\delta]}
\rho(x,t)e^{\sqrt{\tau}x\cdot\omega-\tau t}dxdt
=O(e^{-\tau\delta/2})
\tag {3.5}
$$
as $\tau\longrightarrow\infty$.

Here we prepare

\proclaim{\noindent Lemma 3.1.}
As $\tau\longrightarrow\infty$ we have
$$
\displaystyle
K_1e^{\sqrt{\tau}K_2}\tau^{-(p+1)}
\le e^{\tau T_0}
\vert\int_{D\cap\,(\Bbb R^n\times\,[T_0,\,T_0+\delta])}
\rho(x,t)e^{\sqrt{\tau}x\cdot\omega-\tau t}dxdt\vert
\le K_3e^{\sqrt{\tau}K_4}
\tag {3.6}
$$
where $K_1, K_2, K_3 ,K_4$ are constant and $K_1, K_3>0$.

\endproclaim

{\it\noindent Proof.}

We just describe the key point of the proof of the left half of the inequalities in the case when
$$\displaystyle
\rho(x,t)\ge C_2\,\,\text{for almost all}\,(x,t)\in\,D\cap(\Bbb R^n\times\,[T_0,\,T_0+\delta]).
$$
We have
$$\begin{array}{c}
\displaystyle
e^{\tau T_0}
\int_{D\cap\,(\Bbb R^n\times\,[T_0,\,T_0+\delta])}
\rho(x,t)e^{\sqrt{\tau}x\cdot\omega-\tau t}dxdt
\ge C_2 e^{\tau T_0}
\int_{D\cap\,(\Bbb R^n\times\,[T_0,\,T_0+\delta])}e^{\sqrt{\tau}x\cdot\omega-\tau t}dxdt\\
\\
\displaystyle
=C_2\int_{T_0}^{T_0+\delta}e^{-\tau(t-T_0)}
\left(\int_{D(t-T_0)}e^{\sqrt{\tau}x\cdot\omega}dx\right)dt\\
\\
\displaystyle
=C_2\int_0^{\delta}e^{-\tau s}
\left(\int_{D(s)}e^{\sqrt{\tau}x\cdot\omega}dx\right)ds\\
\\
\displaystyle
\ge C_2e^{\displaystyle
\sqrt{\tau}\inf_{x\in\,\Omega}x\cdot\omega}\int_0^{\delta}\vert D(s)\vert e^{-\tau s}ds\\
\\
\displaystyle
\ge C_1C_2 e^{\displaystyle\sqrt{\tau}\inf_{x\in\Omega}x\cdot\omega}
\int_0^{\delta}s^{p}e^{-\tau s}ds\\
\\
\displaystyle
=C_1C_2\tau^{-(p+1)}e^{\displaystyle
\sqrt{\tau}\inf_{x\in\Omega}x\cdot\omega}\int_0^{\tau\delta}\xi^pe^{-\xi}d\xi.
\end{array}
$$

\noindent
$\Box$

\noindent
Now Theorem 2.1 is a consequence of (3.1) to (3.6).

{\bf\noindent 3.2. A generalization}

In this subsection we give a generalization of Theorem 2.1
to the heat equation with variable coefficients:
$$\displaystyle
au_t=\nabla\cdot\gamma\nabla u+F(x,t)
$$
where $a=a(x)$ is a smooth function with positive values and
$\gamma=\gamma(x)$ is a $n\times n$-real symmetric positive
definite matrix valued smooth function. We apply the idea of the
enclosure method to this equation. For the purpose we construct a
solution with large parameter $\tau$ of the equation
$$
\displaystyle
a v_t+\nabla\cdot\gamma\nabla v=0
\tag {3.7}
$$
which plays the role of the function $e^{\sqrt{\tau}x\cdot\omega-\tau t}$
in Theorem 2.1.

First we consider how to construct a solution of the equation (3.7) in the form:
$$\displaystyle
v(x,t)=e^{-\tau t}w(x).
$$
Substituting this into (3.7), we have
$$
\displaystyle
\nabla\cdot\gamma\nabla w-\tau a w=0.
\tag {3.8}
$$
Note that $\tau>0$ unlike the case when the governing equation is
the wave equation (\cite{ES}).
Let $\varphi$ be a smooth function.
The change of the {\it dependent} variable
$$\displaystyle
w=e^{\sqrt{\tau}\varphi}w^{'}
$$
gives
$$
\displaystyle
\nabla\cdot\gamma\nabla w'
+\sqrt{\tau}\left(2\gamma\nabla\varphi\cdot\nabla w'
+(\nabla\cdot\gamma\nabla\varphi)w'\right)
+\tau(\left(\gamma\nabla\varphi\cdot\nabla\varphi)-a\right)w'=0.
$$
Thus if $\varphi$ satisfies the {\it eikonal equation}
$$\displaystyle
\gamma\nabla\varphi\cdot\nabla\varphi=a,
\tag {3.9}
$$
then $w'$ satisfies
$$
\displaystyle
\nabla\cdot\gamma\nabla w'
+\sqrt{\tau}\left(2\gamma\nabla\varphi\cdot\nabla w'
+(\nabla\cdot\gamma\nabla\varphi)w'\right)=0.
\tag {3.10}
$$
The next lemma is crucial for the construction of the {\it exact}
solution from the asymptotic solution.

\proclaim{\noindent Proposition 3.2.}
Let $\tau$ be an arbitrary positive number.
Given $f\in L^2(\Omega)$ there exists a unique weak solution $u\in H^1(\Omega)$
of the elliptic problem
$$\begin{array}{c}
\displaystyle
\nabla\cdot\gamma\nabla u
+\sqrt{\tau}\left(2\gamma\nabla\varphi\cdot\nabla u
+(\nabla\cdot\gamma\nabla\varphi)u\right)=f\,\,\text{in}\,\Omega,\\
\\
\displaystyle
u=0\,\,\text{on}\,\partial\Omega.
\end{array}
$$
Moreover $u$ has the estimate
$$
\displaystyle
\Vert u\Vert_{H^1(\Omega)}\le C\Vert f\Vert_{L^2(\Omega)}
\tag {3.11}
$$
where $C$ is a positive constant independent of $\tau$.

\endproclaim

{\it\noindent Proof.} Using a change of dependent variable, one
can easily deduce the uniqueness and existence of the solution for
the corresponding fact for the equation (3.8) in $\Omega$. The
positivity of $\tau$ is essential. The problem is the estimate
(3.11). Since
$$\displaystyle
(2\gamma\nabla\varphi\cdot\nabla u)u=\gamma\nabla\varphi\cdot\nabla(u^2),
$$
integration by parts gives
$$\displaystyle
\int_{\Omega}\left(2\gamma\nabla \varphi\cdot\nabla u+(\nabla\cdot\gamma\nabla\varphi)u\right)udx=0.
$$
Then from the equation for $u$ we obtain
$$\displaystyle
\int_{\Omega}\gamma\nabla u\cdot\nabla udx=-\int_{\Omega}fudx.
$$
A standard argument gives (3.11).

\noindent
$\Box$

\noindent
From the equation one knows that the solution $u$ in Proposition 3.2 satisfies
$\nabla\cdot\gamma\nabla u\in L^2(\Omega)$.
A combination of the standard elliptic estimate for the operator $\nabla\cdot\gamma\nabla$ and (3.11)
yields
$$\displaystyle
\Vert u\Vert_{H^2(\Omega)}
\le
C'(\sqrt{\tau}+1)\Vert f\Vert_{L^2(\Omega)}.
$$
Then the Sobolev imbedding yields
$$\displaystyle
\Vert u\Vert_{L^{\infty}(\Omega)}\le C^{''}(\sqrt{\tau}+1)\Vert f\Vert_{L^2(\Omega)}.
\tag {3.12}
$$

Given $N\ge 1$ set
$$\displaystyle
w_N'(x)=\sum_{j=0}^N\frac{a_j(x)}{(\sqrt{\tau})^j}
$$
where $a_0,\cdots,a_N$ satisfy the {\it transport equations}
$$\begin{array}{c}
\displaystyle
2\gamma\nabla\varphi\cdot\nabla a_0
+(\nabla\cdot\gamma\nabla\varphi)a_0=0,\\
\\
\displaystyle
2\gamma\nabla\varphi\cdot\nabla a_{j+1}
+(\nabla\cdot\gamma\nabla\varphi)a_{j+1}=-\nabla\cdot\gamma\nabla a_j,\,\,
j=0,\cdots,N-1.
\end{array}
\tag {3.13}
$$
Then we see that
$$\displaystyle
\nabla\cdot\gamma\nabla w_N'
+\sqrt{\tau}(2\gamma\nabla\varphi\cdot\nabla w_N'
+(\nabla\cdot\gamma\nabla\varphi)w_N')
=\frac{\nabla\cdot\gamma\nabla a_N(x)}{(\sqrt{\tau})^N}.
$$
This means that $w'_N$ is an {\it asymptotic solution} of the equation (3.10).
Then Proposition 3.2 ensures the unique solvability of the problem
$$\begin{array}{c}
\displaystyle
\nabla\cdot\gamma\nabla u
+\sqrt{\tau}\left(2\gamma\nabla\varphi\cdot\nabla u
+(\nabla\cdot\gamma\nabla\varphi)u\right)=
-\frac{\nabla\cdot\gamma\nabla a_N(x)}{(\sqrt{\tau})^N}
\,\,\text{in}\,\Omega,\\
\\
\displaystyle
u=0\,\,\text{on}\,\partial\Omega.
\end{array}
$$
Write $u=R_N$.  Then the function
$$\displaystyle
w^{''}_N=w_N'+R_N
$$
is an exact solution of the equation (2.10) and, by virtue of the estimate (3.11), we have
$$\displaystyle
\Vert R_N\Vert_{H^1(\Omega)}=O\left(\frac{1}{(\sqrt{\tau})^N}\right).
$$
Thus we have
$$\displaystyle
w^{''}_N=\sum_{j=0}^{N-1}\frac{a_j(x)}{(\sqrt{\tau})^j}+O\left(\frac{1}{(\sqrt{\tau})^N}\right).
$$
Note that $w^{''}_N$ is the unique solution of the problem
$$\begin{array}{c}
\displaystyle
\nabla\cdot\gamma\nabla u
+\sqrt{\tau}\left(2\gamma\nabla\varphi\cdot\nabla u
+(\nabla\cdot\gamma\nabla\varphi)u\right)=0\,\,\text{in}\,\Omega,\\
\\
\displaystyle
u=w_N'\,\,\text{on}\,\partial\Omega.
\end{array}
$$

Summing up, we have

\proclaim{\noindent Theorem 3.3.}
Let $\varphi$ satisfy the eikonal equation (3.9).
Given $N\ge 1$ let the functions $a_0,\cdots, a_N$ satisfy
the transport equations (3.13).
Then the unique solution of the elliptic problem
$$\begin{array}{c}
\displaystyle
\nabla\cdot\gamma\nabla w-\tau a w=0\,\,\text{in}\,\Omega,\\
\\
\displaystyle
w=e^{\sqrt{\tau}\varphi}\sum_{j=0}^N\frac{a_j}{(\sqrt{\tau})^j}\,\,\text{on}\,\partial\Omega
\end{array}
$$
satisfies, as $\tau\longrightarrow\infty$
$$\displaystyle
\Vert
e^{-\sqrt{\tau}\varphi}w-\sum_{j=0}^{N-1}\frac{a_j}{(\sqrt{\tau})^j}\Vert_{H^1(\Omega)}
=O\left(\frac{1}{(\sqrt{\tau})^N}\right).
$$
\endproclaim

\noindent
Note that, from (3.12) we also have
$$\displaystyle
\Vert
e^{-\sqrt{\tau}\varphi}w-\sum_{j=0}^{N-2}\frac{a_j}{(\sqrt{\tau})^j}\Vert_{L^{\infty}(\Omega)}
=O\left(\frac{1}{(\sqrt{\tau})^{N-1}}\right)
\tag {3.14}
$$
where $N\ge 2$.

Thus we have a special solution $v$ of the backward heat equation (3.7) in $\Omega\times\Bbb R^+$
having the form
$$\displaystyle
v(x,t)
=e^{-\tau t}e^{\sqrt{\tau}\varphi}
\left\{\sum_{j=0}^{N-1}\frac{a_j(x)}{(\sqrt{\tau})^j}
+O\left(\frac{1}{(\sqrt{\tau})^{N}}\right)\right\}
\tag {3.15}
$$
as $\tau\longrightarrow\infty$ provided $\varphi$, $a_0,\cdots,a_N$
are given.

Now we define another indicator function.

{\bf\noindent Definition 3.1.}
Given $s\in\Bbb R$ define the {\it indicator function}
$I_{\varphi}(\tau;s)$ by the formula
$$\displaystyle
I_{\varphi}(\tau;s)
=e^{\tau s}\int_0^{T}\int_{\partial\Omega}\left(\frac{\partial v}{\partial\nu}u
-\frac{\partial u}{\partial\nu}v\right)dSdt,\,\,\tau>0
$$
where $v$ is given by (3.15).

\noindent

\noindent
For the description of the result we introduce the assumption instead of (A.2):

(A.2)'  There exist positive numbers
$\delta<T-T_0$, $C_1$, $C_2$ and $p\in\,[0,\,\infty[$
such that the $n$-dimensional Lebesgue
measure of the set $D(s)\equiv D\cap(\Bbb R^n\times\{T_0+s\})$
has the lower estimate:
$$\displaystyle
\vert D(s)\vert\ge C_1s^p\,\,\text{for almost all}\, s\in\,[0,\,\delta_1];
$$
the strength of the source $\rho$ multiplied by $a_0(x)$ satisfies
$$\displaystyle
\rho(x,t)a_0(x)\ge C_2\,\,\text{for almost all}\,(x,t)\in\,D\cap(\Bbb R^n\times\,[T_0,\,T_0+\delta])
$$
or
$$\displaystyle
-\rho(x,t)a_0(x)\ge C_2\,\,\text{for almost all}\, (x,t)\in\,D\cap(\Bbb R^n\times\,[T_0,\,T_0+\delta]).
$$

\noindent
We obtain

\proclaim{\noindent Theorem 3.4.}
Assume that the assumptions (A.1), (A.2)' are satisfied.
Then, as $\tau\longrightarrow\infty$ we have
$$\begin{array}{c}
\displaystyle
\lim_{\tau\longrightarrow\infty}\frac{\log\vert I_{\varphi}(\tau;0)\vert}
{\tau}
=-T_0
\end{array}
$$

\noindent
Moreover, we have:

if $s<T_0$, then $\lim_{\tau\longrightarrow\infty}\vert
I_{\varphi}(\tau;s)\vert=0$;

if $s>T_0$, then $\lim_{\tau\longrightarrow\infty}\vert
I_{\varphi}(\tau;s)\vert=\infty$.

\endproclaim

{\it\noindent Proof.}
(3.12) for $N=2$ gives the estimate
$$\displaystyle
\Vert w-a_0(\,\cdot\,)\Vert_{L^{\infty}(\Omega)}=O\left(\frac{1}{\sqrt{\tau}}\right).
$$
Then, it follows from (A.2)' that there exists $\tau_0>0$ such that, for all $\tau\ge\tau_0$
$$\displaystyle
\rho(x,t)w(x)\ge C_2/2\,\,\text{for almost all}\,(x,t)\in\,D\cap(\Bbb R^n\times\,[T_0,\,T_0+\delta])
$$
or
$$\displaystyle
-\rho(x,t)w(x)\ge C_2/2\,\,\text{for almost all}\, (x,t)\in\,D\cap(\Bbb R^n\times\,[T_0,\,T_0+\delta]).
$$

\noindent
Hereafter we take the same course as the proof of Theorem 2.1.

\noindent
$\Box$

{\bf\noindent Example 3.1.}
Consider the simplest case:
$a\equiv 1$ and $\gamma\equiv I_n$.

\noindent Given an arbitrary point $p\in\Bbb
R^n\setminus\overline\Omega$ the function $\varphi$ given by the
formula
$$\displaystyle
\varphi(x;p)=\pm\vert x-p\vert,\,x\in\overline\Bbb R^n\setminus\{p\}
$$
satisfies the equation (3.9) for $a$ and $\gamma$ specified above.
Since
$$\displaystyle
\triangle\varphi=\pm\frac{n-1}{\vert x-p\vert},
$$
the transport equation for $a_0$ in (3.12) becomes
$$\displaystyle
(x-p)\cdot\nabla a_0+\frac{n-1}{2}a_0=0
$$
and the {\it method of bicharacteristics} yields
$$
\displaystyle
a_0(x)
=a_0(\frac{x-p}{\vert x-p\vert}+p)\vert x-p\vert^{-(n-1)/2}.
$$
Thus if we specify the boundary value of $a_0$ on the unit sphere
centered at $p$ as
$$
a_0(x)=1\,\,\text{on}\,\vert x-p\vert=1,
$$
then we have
$$
a_0(x)=\vert x-p\vert^{-(n-1)/2}.
$$
One can also obtain $a_1$ explicitly.
The point is: the function $a_0$ specified above never vanishes on $\overline\Omega$.  From Theorem 3.3 one obtains
the solution $v=v_{\pm}(x,t;p)$ of the backward heat equation having the form
$$\displaystyle
v_{\pm}(x,t;p)=e^{-\tau t}e^{\pm\sqrt{\tau}\vert x-p\vert}
\left\{\vert x-p\vert^{-(n-1)/2}+O\left(\frac{1}{\sqrt{\tau}}\right)\right\}
$$
as $\tau\longrightarrow\infty$.

\section{Proofs of Theorems 2.2 and 2.3}

{\bf\noindent 4.1.  Asymptotic behaviour of an integral}

In this subsection, we consider the case when $n\ge 2$.
we assume that $D\subset\Bbb R^{n+1}$ is a
finite cone with a vertex at $p$ and a bottom face $Q$. More
precisely, $D$ takes the form
$$\displaystyle
D=\cup_{0<s<\delta}\{p+\frac{s}{\delta}(y-p)\,\vert\,y\in Q\}
$$
and $Q$ is a bounded open subset of the $n$-dimensional hyper plane
$$
\left(\begin{array}{c} x\\ t\end{array}\right)\cdot\omega(c)=p\cdot\omega(c)-\delta
$$
where $\delta$ is a positive number.

\noindent
We consider the integral of the function $v(x,t)$ given by (2.2) over $D$:
$$\displaystyle
I(\tau)
=\int_D v(x,t)dxdt
$$

\proclaim{\noindent Lemma 4.1.}
The limit
$$\displaystyle
\lim_{\tau\longrightarrow\infty}\frac{2}{n!}
(c\tau)^{n+1}\exp{\left\{-\sqrt{c^2+1}\tau p\cdot\omega(c)\right\}}
\exp{\left\{-ic\tau\sqrt{1-\frac{1}{c^2\tau}}p\cdot\left(\begin{array}{c}\omega^{\perp}\\
0
\end{array}
\right)\right\}}I(\tau)
$$
exists and has the integral representation:
$$\displaystyle
K_D=2\delta\int_Q
\frac{dS(y)}
{\displaystyle
\left(\frac{\delta\sqrt{c^2+1}}{c}-i(y-p)\cdot\left(\begin{array}{c}
\omega^{\perp}\\ 0
\end{array}
\right)\right)^{n+1}}.
\tag {4.1}
$$

\endproclaim

{\it\noindent Proof.}

Write
$$\begin{array}{c}
\displaystyle
\exp{\left\{-\sqrt{c^2+1}\tau p\cdot\omega(c)\right\}}
\exp{\left\{-ic\tau\sqrt{1-\frac{1}{c^2\tau}}p\cdot\left(\begin{array}{c}
\omega^{\perp}\\ 0\end{array}\right)\right\}}I(\tau)\\
\\
\displaystyle
=\int_D\exp{\sqrt{c^2+1}\tau\left(\left(\begin{array}{c} x\\ t\end{array}\right)-p\right)
\cdot\omega(c)}
\exp{ic\tau\sqrt{1-\frac{1}{c^2\tau}}\left(\left(\begin{array}{c} x\\ t\end{array}\right)-p\right)
\cdot\left(\begin{array}{c} \omega^{\perp}\\ 0\end{array}\right)}dxdt\\
\\
\displaystyle
=\int_0^{\delta}ds\int_Q dS(y)
\left(\frac{s}{\delta}\right)^n
e^{-\sqrt{c^2+1}\tau s}
\exp{\frac{ic\tau s}{\delta}\sqrt{1-\frac{1}{c^2\tau}}(y-p)\cdot\left(\begin{array}{c} \omega^{\perp}
\\ 0\end{array}\right)}\\
\\
\displaystyle
=\delta^{-n}\tau^{-(n+1)}
\int_Q dS(y)
\int_0^{\tau\delta}d\xi
\xi^n e^{-\sqrt{c^2+1}\,\xi}
\exp{\frac{ic\xi}{\delta}
\sqrt{1-\frac{1}{c^2\tau}}(y-p)\cdot\left(\begin{array}{c}
\omega^{\perp}\\
0\end{array}\right)}.
\end{array}
$$

\noindent
Then, a combination of Lebesgue's dominated convergence theorem
and the formula
$$\displaystyle
\int_0^{\infty}\xi^ne^{-\xi}e^{ia\xi}d\xi=\frac{n!}{(1-ia)^{n+1}},\,\,a\in\Bbb R
$$
gives
$$\begin{array}{c}
\displaystyle
\tau^{n+1}
\exp{\left\{-\sqrt{c^2+1}\tau p\cdot\omega(c)\right\}}
\exp{\left\{-ic\tau\sqrt{1-\frac{1}{c^2\tau}}p\cdot\left(\begin{array}{c}
\omega^{\perp}\\ 0\end{array}\right)\right\}}I(\tau)\\
\\
\displaystyle
\longrightarrow
\delta^{-n}
\int_Q dS(y)
\int_0^{\infty}d\xi\,\xi^ne^{-\sqrt{c^2+1}\,\xi}
\exp{\frac{ic\xi}{\delta}(y-p)\cdot\left(\begin{array}{c}
\omega^{\perp}\\ 0\end{array}\right)}\\
\\
\displaystyle
=\delta^{-n}(\sqrt{c^2+1})^{-(n+1)}
\int_Q dS(y)
\left\{\int_0^{\infty}\xi^ne^{-\xi}
\exp{\frac{ic\xi}{\delta\sqrt{c^2+1}}(y-p)\cdot\left(\begin{array}{c}
\omega^{\perp}\\ 0\end{array}\right)}d\xi\right\}\\
\\
\displaystyle
=\frac{n!}{\delta^n(\sqrt{c^2+1})^{n+1}}
\int_Q\frac{dS(y)}
{\displaystyle\left(1-\frac{ic}{\delta\sqrt{c^2+1}}(y-p)\cdot\left(\begin{array}{c}
\omega^{\perp}\\
0
\end{array}
\right)\right)^{n+1}}\\
\\
\displaystyle
=\frac{n!}{2}c^{-(n+1)}K_D.
\end{array}
$$

\noindent
$\Box$

{\bf\noindent 4.2.  Non vanishing of $K_D$}

It seems difficult to check $K_D\not=0$ by using (4.1) directly.
In this subsection, we give an alternative expression of $K_D$ which yields $K_D\not=0$ directly in the case when
$n=2$ provided $D\subset\Bbb R^{2+1}$ is given by the interior of
a tetrahedron with a vertex at $P$ with the bottom face $Q$ located on
the plane $(x,t)^T\cdot\omega(c)=p\cdot\omega(c)-\delta$.
The set $Q$ is given by the interior of the triangle with vertices $y_1,y_2,y_3$
which satisfy $y_j\cdot\omega(c)=p\cdot\omega(c)-\delta$.

Set
$$
\displaystyle
\Delta=\{(\alpha,\beta)\,\vert\,\alpha,\beta\ge 0,\,\alpha+\beta\le 1\}.
$$
Define the maps $\varphi_1,\varphi_2,\varphi_3:\Delta\longrightarrow\Bbb R^{2+1}$
by the formulae
$$\begin{array}{c}
\displaystyle
\varphi_1(\alpha,\beta)=p+\alpha(y_1-p)+\beta(y_2-p),\\
\\
\displaystyle
\varphi_2(\alpha,\beta)=p+\alpha(y_2-p)+\beta(y_3-p),\\
\\
\displaystyle
\varphi_3(\alpha,\beta)=p+\alpha(y_3-p)+\beta(y_1-p)
\end{array}
$$
and set
$$\displaystyle
\Delta_j=\varphi_j(\Delta).
$$
Renumbering $y_1,y_2,y_3$ if necessary, one has the decomposition of $\partial D$:
$$\displaystyle
\partial D=Q\cup\Delta_1\cup\Delta_2\cup\Delta_3.
$$
$Q$ and $\Delta_j,\,j=1,2,3$ satisfy $Q\cup\Delta_j=\emptyset$
and $\text{Int}\,\Delta_j\cap\text{Int}\,\Delta_{j'}=\emptyset$
if $j\not=j'$.

\noindent Let $\nu$ be the unit outward normal vector field to
$\partial D$.  Since $\nu$ takes a constant vector on each
$\Delta_j$, we denote the vector by $\nu_j$.  For simplicity of
description we identify $3+1$ with $1$ and $1-1$ with $3$.
Renumbering $y_1,y_2,y_3$ if necessary, one may assume, in
advance, that
$$\begin{array}{c}
\displaystyle
(\nu_1\times\nu_3)\cdot\omega(c)<0,\\
\\
\displaystyle
(\nu_2\times\nu_1)\cdot\omega(c)<0,\\
\\
\displaystyle
(\nu_3\times\nu_2)\cdot\omega(c)<0.
\end{array}
$$

\proclaim{\noindent Proposition 4.2.}
The formula
$$\displaystyle
K_D\left(\begin{array}{c} c(\omega+i\omega^{\perp})\\
-1
\end{array}
\right)
=c^3\sum_{j=1}^3
\frac{\vert(\nu_j\times\nu_{j-1})\times(\nu_{j+1}\times\nu_j)\vert}
{(\nu_j\times\nu_{j-1})\cdot\vartheta\,
(\nu_{j+1}\times\nu_j)\cdot\vartheta}\nu_j,
\tag {4.2}
$$
is valid where
$$\displaystyle
\vartheta=\left(\begin{array}{c} c(\omega+i\omega^{\perp})\\
-1\end{array}
\right).
$$

\endproclaim

{\it\noindent Proof.} Let $\text{\boldmath $a$}$ be an arbitrary
constant complex vector in four-dimensions. Since
$$\displaystyle
\nabla_{(x,t)}\cdot(v\text{\boldmath $a$})
=\left\{\left(\begin{array}{c} z\\ -\tau\end{array}\right)\cdot\text{\boldmath $a$}\right\}v,
$$
we have
$$\begin{array}{c}
\displaystyle
\left(\begin{array}{c} z\\ -\tau\end{array}\right)\cdot
\text{\boldmath $a$}
\int_D vdxdt
=\sum_{j=1}^3\text{\boldmath $a$}\cdot\nu_j\int_{\Delta_j}vdS(y)
-\text{\boldmath $a$}\cdot\omega(c)\int_Q vdS(y).
\end{array}
\tag {4.3}
$$

\noindent
It is easy to see that, as $\tau\longrightarrow\infty$
$$\displaystyle
\int_Q vdS(y)=O\left(e^{\sqrt{c^2+1}\tau (p\cdot\omega(c)-\delta)}\right).
\tag {4.4}
$$

On the other hand, using a similar computation as in \cite{I5}, we
have
$$\begin{array}{c}
\displaystyle
\int_{\Delta_j}vdS(y)\\
\\
\displaystyle
=\exp{\left\{\sqrt{c^2+1}\tau p\cdot\omega(c)\right\}}
\exp{\left\{ic\tau\sqrt{1-\frac{1}{c^2\tau}}p\cdot\left(\begin{array}{c} \omega^{\perp}\\ 0\end{array}\right)
\right\}}
\frac{\delta^2\vert(\nu_j\times\nu_{j-1})\times(\nu_{j+1}\times\nu_j)\vert}
{(\nu_j\times\nu_{j-1})\cdot\omega(c)\,(\nu_{j+1}\times\nu_j)\cdot\omega(c)}\\
\\
\displaystyle
\times
\int_{\Delta}
\exp{\left\{-\sqrt{c^2+1}\tau(\alpha+\beta)\delta\right\}}
\exp{\left\{-ic\tau\sqrt{1-\frac{1}{c^2\tau}}\delta(\alpha c_1+\beta c_2)\right\}}d\alpha d\beta
\end{array}
\tag {4.5}
$$
where
$$
\displaystyle
c_1=\frac{(\nu_j\times\nu_{j-1})\cdot\left(\begin{array}{c} \omega^{\perp}\\ 0\end{array}\right)}
{(\nu_j\times\nu_{j-1})\cdot\omega(c)},\,\,
c_2=\frac{(\nu_{j+1}\times\nu_j)\cdot\left(\begin{array}{c} \omega^{\perp}\\ 0\end{array}\right)}
{(\nu_{j+1}\times\nu_j)\cdot\omega(c)}.
$$

A combination of the change of variables and Lebesgue's dominated
convergence theorem gives
$$\begin{array}{c}
\displaystyle
(\delta\tau)^2\int_{\Delta}
\exp{\left\{-\sqrt{c^2+1}\tau(\alpha+\beta)\delta\right\}}
\exp{\left\{-ic\tau\sqrt{1-\frac{1}{c^2\tau}}\delta(\alpha c_1+\beta c_2)\right\}}d\alpha d\beta\\
\\
\displaystyle
=\int_0^{\tau\delta}d\beta\int_0^{\tau\delta-\beta}d\alpha
\exp{\left\{-\sqrt{c^2+1}(\alpha+\beta)\right\}}
\exp{\left\{-ic\sqrt{1-\frac{1}{c^2\tau}}(\alpha c_1+\beta c_2)\right\}}\\
\\
\displaystyle
\longrightarrow
\int_0^{\infty}d\beta\int_0^{\infty}d\alpha
\exp{\left\{-\sqrt{c^2+1}(\alpha+\beta)\right\}}
\exp{\left\{-ic(\alpha c_1+\beta c_2)\right\}}\\
\\
\displaystyle
=(c^2+1)^{-1}
\int_0^{\infty}d\beta\int_0^{\infty}d\alpha
\exp{\left\{-(\alpha+\beta)\right\}}
\exp{\left\{-i\frac{c}{\sqrt{c^2+1}}(\alpha c_1+\beta c_2)\right\}}\\
\\
\displaystyle
=(c^2+1)^{-1}
\frac{1}
{\displaystyle\left(1+i\frac{c}{\sqrt{c^2+1}}c_1\right)
\left(1+i\frac{c}{\sqrt{c^2+1}}c_2\right)}\\
\\
\displaystyle
=\frac{1}
{\displaystyle
\left(\sqrt{c^2+1}+icc_1\right)
\left(\sqrt{c^2+1}+icc_2\right)}\\
\\
\displaystyle
=\frac{(\nu_j\times\nu_{j-1})\cdot\omega(c)\,(\nu_{j+1}\times\nu_j)\cdot\omega(c)}
{(\nu_j\times\nu_{j-1})\cdot\vartheta\,(\nu_{j+1}\times\nu_j)\cdot\vartheta}.
\end{array}
\tag {4.6}
$$
A combination of (4.5) and (4.6) gives
$$\begin{array}{c}
\displaystyle
\tau^2\exp{\left\{-\sqrt{c^2+1}\tau p\cdot\omega(c)\right\}}
\exp{\left\{-ic\tau\sqrt{1-\frac{1}{c^2\tau}}p\cdot\left(\begin{array}{c} \omega^{\perp}\\ 0\end{array}\right)
\right\}}\int_{\Delta_j}vdS(y)\\
\\
\displaystyle
\longrightarrow
\frac{\vert(\nu_j\times\nu_{j-1})\times(\nu_{j+1}\times\nu_j)\vert}
{(\nu_j\times\nu_{j-1})\cdot\vartheta\,(\nu_{j+1}\times\nu_j)\cdot\vartheta}.
\end{array}
$$
Then from this, (4.3) and (4.4) we have
$$\begin{array}{c}
\displaystyle
\tau^3\left(\begin{array}{c}
c(\omega+i\sqrt{1-\frac{1}{c^2\tau}}\omega^{\perp})\\
-1\end{array}\right)\cdot\text{\boldmath $a$}
\\
\\
\displaystyle
\times
\exp{\left\{-\sqrt{c^2+1}\tau p\cdot\omega(c)\right\}}
\exp{\left\{-ic\tau\sqrt{1-\frac{1}{c^2\tau}}p\cdot\left(\begin{array}{c}\omega^{\perp}\\ 0\end{array}\right)
\right\}}I(\tau)\\
\\
\displaystyle
\longrightarrow
\sum_{j=1}^3
\frac{\vert(\nu_j\times\nu_{j-1})\times(\nu_{j+1}\times\nu_j)\vert\nu_j\cdot\text{\boldmath $a$}}
{(\nu_j\times\nu_{j-1})\cdot\vartheta\,
(\nu_{j+1}\times\nu_j)\cdot\vartheta}.
\end{array}
\tag {4.7}
$$

A combination of Lemma 4.1 and (4.7) yields
$$\begin{array}{c}
\displaystyle
K_D
\left(\begin{array}{c} c(\omega+i\omega^{\perp})\\
-1\end{array}
\right)\cdot\text{\boldmath $a$}\\
\\
\displaystyle
=\lim_{\tau\longrightarrow\infty}
\left(\begin{array}{c}
c(\omega+i\sqrt{1-\frac{1}{c^2\tau}}\omega^{\perp})\\
-1\end{array}
\right)\cdot\text{\boldmath $a$}\\
\\
\displaystyle
\times
\frac{2}{2!}(c\tau)^3\exp{\left\{-\sqrt{c^2+1}\tau p\cdot\omega(c)\right\}}
\exp{\left\{-ic\sqrt{1-\frac{1}{c^2\tau}}\tau
p\cdot\left(\begin{array}{c} \omega^{\perp}\\ 0\end{array}\right)\right\}}
I(\tau)\\
\\
\displaystyle
=
c^3\sum_{j=1}^3
\frac{\vert(\nu_j\times\nu_{j-1})\times(\nu_{j+1}\times\nu_j)\vert\nu_j\cdot\text{\boldmath $a$}}
{(\nu_j\times\nu_{j-1})\cdot\vartheta\,
(\nu_{j+1}\times\nu_j)\cdot\vartheta}.
\end{array}
$$
Since $\text{\boldmath $a$}$ is arbitrary, we obtain (4.2).

\noindent
$\Box$

\noindent
Since $\nu_1,\nu_2,\nu_3$ are linearly independent, one has the conclusion

\proclaim{\noindent Corollary 4.3.}
$K_D\not=0$ for all $c>0$.

\endproclaim

{\bf\noindent 4.3.  Proof of Theorem 2.2}

Integration by parts gives
$$\displaystyle
I_{\omega,\,\omega^{\perp},c}(\tau;0)
=\int_0^T\int_{\Omega}f(x,t)v(x,t)dxdt
-\int_{\Omega}u(x,T)v(x,T)dx.
$$
Then (4.8) yields that the integral
$$
e^{-\tau\sqrt{c^2+1}\,h_D(\omega(c))}\int_{\Omega}u(x,T)v(x,T)dx
$$
is exponentially decaying as $\tau\longrightarrow\infty$.  Thus the integral
$\displaystyle e^{-\tau\sqrt{c^2+1}\,h_D(\omega(c))}I_{\omega,\,\omega^{\perp},c}(\tau;0)$ modulo exponentially
decaying as $\tau\longrightarrow\infty$,
coincides with the integral
$$\begin{array}{c}
\displaystyle
e^{-\tau\sqrt{c^2+1}\,h_D(\omega(c))}\int_0^T\int_{\Omega}f(x,t)v(x,t)dx\\
\\
\displaystyle
=e^{-\tau\sqrt{c^2+1}\,h_D(\omega(c))}
\sum_{j=1}^N\int_{P_j\times]T_j,\,T[}\rho_j(x,t)e^{-\tau t}e^{x\cdot z}dxdt.
\end{array}
\tag {4.9}
$$
Since $\omega(c)$ is regular with respect to $D$, there exits a
unique point $p$ on $\partial D$ such that
$p\cdot\omega(c)=h_D(\omega(c))$.  This $p$ should belong to some
$\displaystyle P_{j_0}\times\{T_{j_0}\}$ since the time component
of $\omega(c)$ is negative.  Then $P_{j'}\times\,]T_{j'},T[$ with
$j'\not=j_0$ should be contained in the half space
$(x,t)^T\cdot\omega(c)<h_D(\omega)-\delta$ with a positive number
$\delta$ independent of $j'$.  Thus the right hand side of (4.9)
modulo exponentially decaying as $\tau\longrightarrow\infty$, coincides with
$$\displaystyle
e^{-\tau\sqrt{c^2+1}\,h_D(\omega(c))}
\int_{D'}\rho_{j_0}(x,t)e^{-\tau t}e^{x\cdot z}dxdt
\tag {4.10}
$$
where
$$\displaystyle
D'=\left(P_{j_0}\times\,]T_j,\,T[\right)\cap\{(x,t)\,\vert\,(x,t)^T\cdot\omega\ge h_D(\omega)-\delta\}.
$$
Choosing a smaller $\delta$ if necessary, one may assume that this
$D'$ is given by the interior of a tetrahedron with vertex $p$ and
has the factorization as that of $D$ in subsection 4.1:
$$\displaystyle
D'=\cup_{0<s<\delta}\{p+\frac{s}{\delta}(y-p)\,\vert\,y\in Q\}
$$
where $Q$ is the interior of a triangle lying on the plane
$(x,t)^T\cdot\omega(c)=h_D(\omega)-\delta$.
This is the key point of the proof.
Divide (4.10) into two parts:
$$\begin{array}{c}
\displaystyle
e^{-\tau\sqrt{c^2+1}\,h_D(\omega(c))}\rho_{j_0}(p)
\int_{D'}e^{-\tau t}e^{x\cdot z}dxdt
+e^{-\tau\sqrt{c^2+1}\,h_D(\omega(c))}
\int_{D'}\{\rho_{j_0}(x,t)-\rho_{j_0}(p)\}e^{-\tau t}e^{x\cdot z}dxdt\\
\\
\displaystyle
\equiv I+II.
\end{array}
$$
From a combination of Lemma 4.1 and Corollary 4.3 one knows that $I$ as $\tau\longrightarrow\infty$
decays really algebraically:
$$\displaystyle
\lim_{\tau\longrightarrow\infty}(c\tau)^3\vert I\vert=\vert K_{D'}\vert\vert\rho_{j_0}(p)\vert.
$$
Using the H\"older continuity of $\rho_{j_0}$, one can also show easily
$$\displaystyle
\tau^3\vert II\vert=O\left(\tau^{-\theta_{j_0}}\right).
$$

\noindent
Summing up, one concludes that, as $\longrightarrow\infty$
$$\displaystyle
(c\tau)^3e^{-\tau\sqrt{c^2+1}\,h_D(\omega(c))}\vert I_{\omega,\,\omega^{\perp},c}(\tau;0)\vert
\longrightarrow \vert K_{D'}\vert\vert\rho_{j_0}(p)\vert.
$$
Then all of the statements follow from this asymptotic formula.
This completes the proof of Theorem 2.2.

{\bf\noindent 4.4.  Proof of Theorem 2.3}

This is the case when $n=1$.  Since
$$\displaystyle
z\cdot z
=\tau+i2c^2\tau^2\sqrt{1-\frac{1}{c^2\tau}},
$$
we have
$$\displaystyle
v(x,t)
=\exp\left\{\sqrt{c^2+1}\tau\left(\begin{array}{c} x\\
t\end{array}\right)\cdot\omega(c)\right\}
\exp\left\{
ic\tau\left(\begin{array}{c} x\\ t
\end{array}\right)\cdot\left(\begin{array}{c} 1\\
-2c\tau\end{array}\right)
\sqrt{1-\frac{1}{c^2\tau}}\right\}.
$$

\noindent Note that this function is oscillatory higher than the
case when $n\ge 2$ because of the existence of the growing factor
$\tau^2$ in the imaginary part of the phase function.

As same as the proof of Theorem 2.2 it suffices to show that the integral
$$\displaystyle
\exp\left\{-\sqrt{c^2+1}\tau h_D(\omega(c))\right\}\int_D v(x,t)\rho(x,t)dxdt
$$
is {\it really algebraic decaying} as $\tau\longrightarrow\infty$.

\noindent
Remarks are in order.

$\bullet$  One may assume that $D$ is given by the interior
of a {\it triangle} with vertices $p$, $y_0$ and $y_1$ with
$p\cdot\omega=h_D(\omega(c))$, $y_0\cdot\omega(c)=y_1\cdot\omega(c)=h_D(\omega(c))-\delta$
and has the form
$$\displaystyle
D=\cup_{0<s<\delta}\{p+\frac{s}{\delta}(y-p)\,\vert\,y\in Q\}
$$
where $Q$ is the segment with endpoints $y_0$, $y_1$.  And also one may assume that $\rho_j, j=,\cdots,N$
are given by a single $\rho\equiv\rho_{j_0}$ for some $j_0$.

$\bullet$  One may assume that
the vectors $y_1-p$ and $y_0-p$ coincide with $\vert y_1-p\vert (0\,\,1)^T$
and $\vert y_0-p\vert (-\text{sgn}\,c\,\,\,0)^T$, respectively.

$\bullet$  Since we have assumed that $\rho$ coincides with a function
of class $C^2$ in a neighbourhood of the point $p$ in $\Bbb R^{1+1}$,
one may assume that, for all $(x,t)\in D$
$$
\displaystyle
\rho(x,t)=\rho(p)+\nabla\rho(p)\cdot\{(x\,\,t)^T-p\}+O(\vert (x\,\,t)^T-p\vert^2)
$$
where $\rho(p)(\not=0)$ and $\nabla\rho(p)$ are the corresponding
values of the $C^2$-extension of the original $\rho$.

\noindent
Then one has
the expression
$$\begin{array}{c}
\displaystyle
\exp\left\{-\sqrt{c^2+1}\tau p\cdot\omega(c)\right\}
\exp\left\{-ic\tau p\cdot\left(\begin{array}{c}
1\\ -2c\tau\end{array}\right)\sqrt{1-\frac{1}{c^2\tau}}\right\}
\int_D v(x,t)\rho(x,t)dxdt\\
\\
\displaystyle
=I(\tau)+II(\tau)+III(\tau)
\end{array}
$$
where
$$\begin{array}{c}
\displaystyle
\frac{I(\tau)}{\rho(p)}
\\
\\
\displaystyle
=\int_D
\exp\left\{\sqrt{c^2+1}\tau
\{\left(\begin{array}{c} x\\
t\end{array}\right)- p\}\cdot\omega(c)\right\}
\exp\left\{ic\tau
\{\left(\begin{array}{c} x\\ t
\end{array}\right)-
p\}\cdot\left(\begin{array}{c}
1\\ -2c\tau\end{array}\right)\sqrt{1-\frac{1}{c^2\tau}}\right\}dxdt,
\end{array}
$$
$$
\begin{array}{c}
\displaystyle
II(\tau)\\
\\
\displaystyle
=\int_D
\exp\left\{\sqrt{c^2+1}\tau
\{\left(\begin{array}{c} x\\
t\end{array}\right)- p\}\cdot\omega(c)\right\}
\exp\left\{ic\tau
\{\left(\begin{array}{c} x\\ t
\end{array}\right)-
p\}\cdot\left(\begin{array}{c}
1\\ -2c\tau\end{array}\right)\sqrt{1-\frac{1}{c^2\tau}}\right\}\\
\\
\times
\nabla\rho(p)\cdot\{(x\,\,t)^T-p\}
dxdt,
\end{array}
$$
and
$$
\begin{array}{c}
\displaystyle
III(\tau)
\\
\\
\displaystyle
=\int_D
\exp\left\{\sqrt{c^2+1}\tau
\{\left(\begin{array}{c} x\\
t\end{array}\right)- p\}\cdot\omega(c)\right\}
\exp\left\{ic\tau
\{\left(\begin{array}{c} x\\ t
\end{array}\right)-
p\}\cdot\left(\begin{array}{c}
1\\ -2c\tau\end{array}\right)\sqrt{1-\frac{1}{c^2\tau}}\right\}\\
\\
\displaystyle
\times O(\vert(x\,\,t)^T-p\vert^2)dxdt.
\end{array}
$$

\noindent
Here we prove that, as $\tau\longrightarrow\infty$
$I(\tau)\sim C\tau^{-3}$ with $C\not=0$ and $II(\tau)=III(\tau)=O(\tau^{-4})$.

First we study the asymptotic behaviour of $I(\tau)$ as $\tau\longrightarrow\infty$.
$$\begin{array}{c}
\displaystyle
I(\tau)/\rho(p)
=\int_0^{\delta}ds
\int_Q
\left(\frac{s}{\delta}\right)
e^{-s\sqrt{c^2+1}\tau}
\exp\left\{
ic\tau\frac{s}{\delta}(y-p)\cdot\left(\begin{array}{c}
1\\ -2c\tau\end{array}\right)\sqrt{1-\frac{1}{c^2\tau}}\right\}dS(y)
\\
\\
\displaystyle
=\frac{1}{\delta\tau^2}
\int_Q dS(y)
\int_0^{\delta\tau}\xi e^{-\xi\sqrt{c^2+1}}
\exp\left\{i\frac{c}{\delta}\xi(y-p)\cdot
\left(\begin{array}{c}
1\\ -2c\tau\end{array}\right)
\sqrt{1-\frac{1}{c^2\tau}}\right\}d\xi.
\end{array}
$$

\noindent
Set
$$
\displaystyle
B(y,\tau)
=\frac{c}{\delta}
(y-p)\cdot\left(\begin{array}{c}
1\\ -2c\tau\end{array}\right)\sqrt{1-\frac{1}{c^2\tau}}.
$$

Since
$$\begin{array}{l}
\displaystyle
\int_0^{\delta\tau}\xi e^{-\xi\sqrt{c^2+1}}
\exp\left\{i\frac{c}{\delta}\xi(y-p)\cdot
\left(\begin{array}{c}
1\\ -2c\tau\end{array}\right)
\sqrt{1-\frac{1}{c^2\tau}}\right\}d\xi\\
\\
\displaystyle
=
\frac{1}{\displaystyle\left( \sqrt{c^2+1}
-iB(y,\tau)\right)^2}
+O(e^{-\tau\delta\sqrt{c^2+1}}),
\end{array}
$$
we obtain
$$
\displaystyle
\tau^2\delta I(\tau)/\rho(p)
=\int_Q
\frac{dS(y)}{\displaystyle\left( \sqrt{c^2+1}
-iB(y,\tau)
\right)^2}
+O(e^{-\tau\delta\sqrt{c^2+1}}).
\tag {4.11}
$$

\noindent Using a parameterization of the segment $Q$, we have

\proclaim{\noindent Lemma 4.4.}
The formulae
$$
\displaystyle
\int_Q
\frac{dS(y)}{\displaystyle\left( \sqrt{c^2+1}
-iB(y,\tau)
\right)^2}
=
\frac{\vert y_1-y_0\vert}
{\displaystyle\left(\sqrt{c^2+1}
-iB(y_1,\tau)\right)\left(\sqrt{c^2+1}
-iB(y_0,\tau)\right)}
\tag {4.12}
$$
and
$$\begin{array}{l}
\displaystyle
\int_Q
\frac{2\nabla\rho(p)\cdot(y-p)}
{\displaystyle\left(\sqrt{c^2+1}
-iB(y,\tau)
\right)^3}dS(y)
=-i\frac{\vert y_1-y_0\vert}
{\displaystyle B(y_1,\tau)-B(y_0,\tau)}
\\
\\
\displaystyle
\times
(
\frac{\nabla\rho(p)\cdot(y_1-p)}
{\displaystyle\left( \sqrt{c^2+1}
-iB(y_1,\tau)\right)^2}
-\frac{\nabla\rho(p)\cdot(y_0-p)}
{\displaystyle\left( \sqrt{c^2+1}
-iB(y_0,\tau)\right)^2}
\\
\\
\displaystyle
-\frac{\nabla\rho(p)\cdot(y_1-y_0)}
{\displaystyle
\left(\sqrt{c^2+1}
-iB(y_1,\tau)\right)
\left(\sqrt{c^2+1}
-iB(y_0,\tau)\right)}),
\end{array}
\tag {4.13}
$$
are valid.
\endproclaim
{\it\noindent Proof.} Let $y=y(\eta)$ be the parameterization of
$Q$ with $y(0)=y_0$, $y(1)=y_1$ and $y'(\eta)=y_1-y_0$. Since
$$\displaystyle
\frac{d}{d\eta}B(y,\tau)=B(y_1,\tau)-B(y_0,\tau),
$$
we have
$$\displaystyle
\frac{d}{d\eta}\left(\frac{1}{\sqrt{c^2+1}-iB(y,\tau)}\right)
=i\frac{B(y_1,\tau)-B(y_2,\tau)}{\left(\sqrt{c^2+1}-iB(y,\tau)\right)^2}
$$
This gives (4.12) and a combination of a similar identity and integration by parts yields also (4.13).

\noindent
$\Box$

\noindent
Since
$$
\displaystyle
B(y_1,\tau)=-\frac{2c^2\vert y_1-p\vert}{\delta}\tau\sqrt{1-\frac{1}{c^2\tau}}
$$
and
$$
\displaystyle
B(y_0,\tau)=-\frac{\vert c\vert\vert y_0-p\vert}{\delta}\sqrt{1-\frac{1}{c^2\tau}},
$$
from (4.12)
we have, as $\tau\longrightarrow\infty$
$$
\displaystyle
\tau\int_Q
\frac{dS(y)}{\displaystyle\left( \sqrt{c^2+1}
-iB(y,\tau)
\right)^2}
\longrightarrow
-i\frac{\delta}{2c^2}
\frac{\vert y_1-y_0\vert^2}
{\displaystyle\vert y_1-p\vert\left(\sqrt{c^2+1}+i\frac{\vert c\vert}{\delta}\vert y_0-p\vert\right)}.
$$
A combination of this and (4.11) gives the formula:
$$
\displaystyle
\lim_{\tau\longrightarrow\infty}
2(c\tau)^2\tau I(\tau)
=
-\frac{i\vert y_1-y_0\vert^2\rho(p)}
{\displaystyle\vert y_1-p\vert\left(\sqrt{c^2+1}+i\frac{\vert c\vert}{\delta}\vert y_0-p\vert\right)}
(\not=0).
$$

For the estimation of $II(\tau)$ one has to make use of the
growing factor $\tau^2$ in the imaginary part of the phase
function in the integrand. For the purpose write
$$\displaystyle
II(\tau)
=\frac{1}{\delta^2\tau^{3}}
\int_Q\nabla\rho(p)\cdot(y-p) dS(y)
\int_0^{\tau\delta}
\xi^2 e^{-\xi\sqrt{c^2+1}}
e^{i\xi B(y,\tau)}d\xi.
$$
Since
$$\displaystyle
\int_0^{\tau\delta}\xi^2 e^{-\xi\sqrt{c^2+1}}
e^{i\xi B(y,\tau)}d\xi
=\frac{2}
{\displaystyle
\left(\sqrt{c^2+1}-iB(y,\tau)\right)^3}
+O((\tau\delta)^2e^{-\tau\delta\sqrt{c^2+1}}),
$$
we have
$$
\displaystyle
\delta^2\tau^3
II(\tau)
=\int_Q
\frac{2\nabla\rho(p)\cdot(y-p)dS(y)}
{\displaystyle\left(\sqrt{c^2+1}
-iB(y,\tau)
\right)^3}
+O((\tau\delta)^2 e^{-\tau\delta\sqrt{c^2+1}}).
$$
Here, from (4.13) we have, as $\tau\longrightarrow\infty$
$$
\displaystyle
\tau\int_Q
\frac{2\nabla\rho(p)\cdot(y-p)dS(y)}
{\displaystyle\left(\sqrt{c^2+1}
-iB(y,\tau)
\right)^3}
\longrightarrow
-\frac{i\delta}{2c^2}\frac{\vert y_1-y_0\vert}{\vert y_1-p\vert}
\frac{\nabla\rho(p)\cdot(y_0-p)}
{\displaystyle
\left(\sqrt{c^2+1}+i\frac{\vert c\vert}{\delta}\vert y_0-p\vert\right)^2}
$$
This gives the estimate $II(\tau)=O(\tau^{-4})$ as $\tau\longrightarrow\infty$.
The estimation of $III(\tau)$ is rather easier than that of $II(\tau)$:
$$\begin{array}{l}
\displaystyle
\vert III(\tau)\vert
\le C\int_0^{\delta}ds\int_Q(\frac{s}{\delta})
e^{-s\tau\sqrt{c^2+1}}
(\frac{s}{\delta})^2\vert y-p\vert^2 dS(y)\\
\\
\displaystyle
\le\frac{C}{\tau^4}\int_{0}^{\infty}(\frac{\xi}{\delta})^3e^{-\xi\sqrt{c^2+1}}d\xi
\int_Q\vert y-p\vert^2dS(y)=O(\tau^{-4}).
\end{array}
$$

Summing up, we conclude the existence of the {\it nonzero limit} of the integral
$$\displaystyle
\tau^3\exp\left\{-\sqrt{c^2+1}\tau p\cdot\omega(c)\right\}
\exp\left\{-ic\tau p\cdot\left(\begin{array}{c}
1\\ -2c\tau\end{array}\right)\sqrt{1-\frac{1}{c^2\tau}}\right\}
\int_D v(x,t)\rho(x,t)dxdt
$$
as $\tau\longrightarrow\infty$.  This completes the proof of Theorem 2.3.

$$\quad$$

\centerline{{\bf Acknowledgement}}

This research was partially supported by Grant-in-Aid for
Scientific Research (C)(2) (No.  15540154) of Japan  Society for
the Promotion of Science.

$$\quad$$

\vskip1cm
\noindent
e-mail address

ikehata@math.sci.gunma-u.ac.jp

\begin{thebibliography}{99}





\bibitem{ABR}   Avdonin, S. A., Belishev, M. I. and Rozhkov, Yu. S.,
              \newblock The BC-method in the inverse problem for the heat equation,
              \newblock J. Inv. Ill-Posed Problems, 5(1997), 309-322.






\bibitem{EH1}  EL Badia, A. and Ha Duong, T.,
              \newblock An inverse source problem in potential analysis,
              \newblock Inverse Problems, 16(2000), 651-663.





\bibitem{EH2} EL Badia, A. and Ha Duong, T.,
             \newblock On an inverse source problem for the heat eqution.  Application to a pollution detection
             problem,
             \newblock J. Inv. Ill-Posed Problems, 10(2002), 585-599.







\bibitem{CD} Cannon, J. R. and DuChateau, P.,
             \newblock Structual identification of an unknown source term in a heat equation,
             \newblock Inverse Problems, 14(1998), 535-551.







\bibitem{CE}  Cannon, J. R. and Esteva, S. P.,
              \newblock An inverse problem for the heat equation,
              \newblock Inverse Problems, 2(1986), 395-403.




\bibitem{ES} Egorov, V., Yu. and Shubin, M. A.,
             \newblock I. Linear Partial Differential Equations.  Elements of the Modern Theory,
             \newblock Partial Differential Equations II, Egorov, V. Yu. and Subin, M. A. (Eds.),
             Encyclopaedia of Mathematical Sciences, 31, 82-86, Springer, 1988.




\bibitem{FI} Fursikov and Imanuvilov,
             \newblock Contorollability of Evolution Equations,
             Lecture Notes Series 34. RIM-GARC, Seoul National
             University, 1996.



\bibitem{GL} Gurariy, V. I. and Lusky, W.,
             \newblock Geometry of M\"untz spaces and related questions,
             Lecture Notes in Math., 1870, Springer, 2005.





\bibitem{HR} Hettlich, F. and Rundell, W.,
             \newblock Identification of a discontinuous source in the heat equation,
             \newblock Inverse Problems, 17(2001), 1465-1482.




\bibitem{I1} Ikehata, M.,
             \newblock Reconstruction of a source domain from the Cauchy data,
             \newblock Inverse Problems, 15(1999), 637-645.



\bibitem{I2} Ikehata, M.,
             \newblock Enclosing a polygonal cavity in a two-dimensional bounded domain from Cauchy data,
             \newblock Inverse Problems, 15(1999), 1231-1241.



\bibitem{I3} Ikehata, M.,
             \newblock On reconstruction in the inverse conductiviy problem with one measurement,
             \newblock Inverse Problems, 16(2000), 785-793.


\bibitem{I4} Ikehata, M.,
            \newblock Reconstruction of the support function for inclusion from boundary measurements,
             \newblock J. Inv. Ill-Posed Problems, 8(2000), 367-378.





\bibitem{I5} Ikehata, M.,
             \newblock Exponentially growing solutions and the Cauchy problem,
             \newblock Appl. Anal., 78(2001), no. 1-2, 79-95.




\bibitem{I6} Ikehata, M.,
             \newblock Inverse scattering problems and the enclosure method,
             \newblock Inverse Problems, 20(2004), 533-551.




\bibitem{I7} Ikehata, M.,
             \newblock An inverse transmission scattering problem and the enclosure method,
             \newblock Computing, 75(2005), 133-156.



\bibitem{I8} Ikehata, M.,
              \newblock The Herglotz wave function, the Vekua transform and the enclosure method,
              \newblock Hiroshima Math. J., 35(2005), 485-506.


\bibitem{IO} Ikehata, M. and Ohe, T.,
             \newblock A numerical method for finding the convex hull of polygonal cavities
             using the enclosure method,
             \newblock Inverse Problems, 18(2002), 111-124.





\bibitem{Is} Isakov, V., Inverse source problems,
             \newblock Mathematical surveys and monographs, No. 34,
             \newblock Americal Mathematical Society, Providence, Rhode Island, 1990.



\bibitem{LR} Lebeau, G and Robbiano, L.,
             Contr\^ole exact de l'\'equation de la chaleur,
             \newblock Comm. PDE, 20(1995), 335-356.




\bibitem{TLA} Trong, D. D., Long, T. N. and Alain, D. N. P.,
             \newblock Nonhomogeneous heat equation: Identification and regularization for the inhomogeneous term,
             \newblock J. Math. Anal. Appl., 312(2005), 93-104.



\bibitem{Y} Yamamoto, M.,
            \newblock Stability, reconstruction formula and regularization for
            an inverse source hyperbolic problem by a control method,
            \newblock Inverse Problems, 11(1995), 481-496.



\bibitem{YO} Yamatani, K. and Ohnaka, K.,
            \newblock An estimation method for point sources of multidimensional diffusion equation,
            \newblock Appl. Math. Modelling, 21(1997), 77-84.











\end{thebibliography}
\end{document}